\documentclass[11pt]{amsart}
\usepackage{amsfonts}

\newtheorem{theorem}{\bf Theorem}

\newtheorem{proposition}{Proposition}[section]
\newtheorem{lemma}{Lemma}

\newtheorem{definition}{Definition}

\usepackage{enumerate,mathrsfs}
\usepackage{amssymb,amsmath,graphicx,amsthm,mathtools,hyperref}

\usepackage{booktabs}

\theoremstyle{plain}
\usepackage{siunitx}

\usepackage{algorithm}
\usepackage{algorithmic,colortbl}

\makeindex
\usepackage[intoc]{nomencl} 

\makenomenclature

\usepackage{hyperref}
\hypersetup{nesting=true,debug=true,naturalnames=true}
\usepackage{graphicx,amssymb,upref}

\usepackage{subcaption}
\usepackage{enumerate,float}
\usepackage{fullpage}

\begin{document}

\title{Tempered Fractional Hawkes Process and Its Generalization}
\author[]{Neha Gupta}
\address{\emph{Department of Mathematics and Statistics,
		Indian Institute of Technology Kanpur, Kanpur 208016, India.}}
\email{nehagpt@iitk.ac.in}
\author[]{Aditya Maheshwari}
\address{\emph{Operations Management and Quantitative Techniques Area, Indian Institute of Management Indore, Indore 453556, India.}}
\email{adityam@iimidr.ac.in}



   \begin{abstract}
   Hawkes process (HP) is a point process with a conditionally dependent intensity function.  This paper defines the tempered fractional Hawkes process (TFHP) by time-changing the HP with an inverse tempered stable subordinator. We obtained results that generalize the fractional Hawkes process defined in \cite{hainaut2020fractional} to a tempered version which has \textit{semi-heavy tailed} decay.
   We derive  the mean, the variance, covariance and the governing fractional difference-differential equations of the TFHP. Additionally, we introduce  the generalized fractional Hawkes process (GFHP) by time-changing the HP with the inverse L\'evy subordinator. This definition encompasses all potential (inverse L\'evy) time changes as specific instances. We also explore the distributional characteristics and the governing difference-differential equation of the one-dimensional distribution for the GFHP.   
\end{abstract}

   \maketitle

\section{Introduction}

Hawkes process (HP) is a point process that describes the occurrences of events over time, where events can trigger further events (affect the arrival rate) in a self-exciting manner. In other words, the intensity rate depends on past event history and can be modelled using a stochastic intensity function with a decaying effect of past events. 
Hawkes process was first introduced in the early seventies by A. G. Hawkes (see \cite{hawkes1971spectra, hawkes1971point}). This process found applications in several areas, for example, in modelling terrorist activities \cite{Porter_2012}, finance (see 
\cite{Bacry_2015,hainaut2017clustered,Hainaut_2019,hautsch2006modelling, bowsher2007modelling,cartea2014buy,chavez2005estimating}), and
seismology (see \cite{hawkes1973cluster,ogata1988statistical,ogata1998space}). For an extensive exploration of the applications of the Hawkes process, please see \cite{reinhart2018review}.\\\\
The decay of intensity function (or influence of past occurrences) is exponential in nature of a common Hawkes process (HP) and assures Markov property  for its intensity function. With the  Markov property, we gain access to a thoroughly explored  toolbox of stochastic calculus, revealing numerous significant findings related to the HP.  
Exponential decay inadequately represents processes with enduring historical memory, whereas a power-law decay model aptly captures long-memory phenomena. Recently, Hainaut (2020) \cite{hainaut2020fractional} defined the fractional Hawkes process by time-changing the HP by an inverse $\beta$-stable subordinator. This process provides an alternative way to model  power-law decay and non-Markov covariance, including long memory property. In this paper, we use tempered $\beta$-stable subordinator to further expand the scope of \cite{hainaut2020fractional}.\\\\
Tempered stable subordinator (TSS) was first studied by \cite{Rosinski09}, and later tempered fractional fractional derivatives were developed by  \cite{temfraccal}. The main idea of exponential tempering is to model \textit{semi-heavy tails}, which behaves like power law at moderate time scales while converging to a Gaussian at long time scales. It offers an intermediate model between the \textit{thin-tail} and \textit{heavy-tail} model. Analogously, the TSS is useful in modelling \textit{semi-long range dependence}, where correlation decays like power law in a moderate time scale and becomes short-range dependence at long time scales. Another advantage of the TSS is that it has finite moments, unlike the $\beta$-stable subordinator. Due to above-mentioned reasons, researchers have studied the TSS extensively in the recent past; for example, the inverse of TSS is studied by \cite{Kumar-ITSS}, the multivariate version is studied by \cite{multi-tss}, the mixture of TSS by \cite{gupta2021stochastic}, tempered negative binomial by \cite{tsfnbp}, and \textit{etc.} \\

In this paper, we define a tempered fractional Hawkes process (TFHP) by time-changing the HP with the first-hitting time or inverse of the TSS. This definition generalizes the fractional Hawkes process (see \cite{hainaut2020fractional}), and  if we put the tempering parameter in the TFHP to zero, we get the fractional Hawkes process. We compute the governing difference-differential equation for the pdf and its Laplace transform of the TFHP. We derive the mean, variance and covariance structure of the TFHP. We also define the generalized fractional Hawkes process (GFHP) using a general L\'evy subordinator as a time-change in HP. Using the general L\'evy subordinator, all possible time-changes become a particular case of the GFHP. We work out the distributional properties of the GFHP  and the governing difference-differential equation for its one-dimensional distributions. 

The rest of the paper is organized as follows. In Section \ref{sec:pre}, we present some results and definitions that are used in later sections. Section \ref{sec:TFHP} defines the TFHP and derives mean, variance and  the corresponding difference-differential equation. In  Section \ref{sec:covstruc}, we find out the covariance structure of the TFHP. Finally in Section \ref{sec:GFHP}, we define and derive results related to the  GFHP.

\section{Preliminaries}\label{sec:pre}
In this section, we define and present some preliminary results and definitions that will be used later in this paper.

\subsection{L\'evy subordinator and its inverse}
A L\'evy subordinator (hereafter referred to as the subordinator) $\{D_{f}(t)\}_{t\geq0}$ is a non-decreasing L\'evy process and its Laplace transform (LT) (see \cite[Section 1.3.2]{appm}) has the form
			\begin{equation*}
				\mathbb{E}[e^{-s D_{f}(t)}]=e^{-tf(s)},
				\;{\rm where}\; 
				f(s)=a+b s+\int_{0}^{\infty}(1-e^{-s x})\upsilon(dx),~b\geq0, s>0,
			\end{equation*}
			is the Bernst\'ein function (see \cite{Bernstein-book} for more details). The above representation is called the L\'evy-Khintchine representation, and ($a$, $b$, $\upsilon$) is the L\'evy triplet of $f$.
			Here, $a, \; b \geq 0$, and $\upsilon(dx)$ is a non-negative L\'evy measure on a positive half-line satisfying 
			\begin{equation*}
				\int_{0}^{\infty}(x\wedge 1)\upsilon(dx)<\infty~~{\rm and}~~\upsilon([0,\infty))= \infty
			\end{equation*}
			which ensures that the sample paths of $\{D_{f}(t)\}_{t \geq 0}$ are almost surely $(a.s.)$  strictly increasing.
			Also, the first-exit time of $\{D_f(t)\}_{t\geq0}$ is defined as
			\begin{equation}\label{inverse-sub}
				E_{f}(t)=\inf\{r\geq 0:D_{f}(r)>t\},
			\end{equation}
			which is the right-continuous inverse of the subordinator $\{D_f(t)\}_{t\geq 0}$.
			The process $\{E_{f}(t)\}_{t \geq 0}$ is non-decreasing, and its sample paths are continuous.\\

  As a special case of the L\'evy subordinator, we consider $$f(s)=(s+\nu)^\beta-\nu^\beta,~s>0,\nu>0,0<\beta<1,$$  it is called as the tempered $\beta$-stable subordinator (TSS). When $\nu=0$, it is called as the $\beta$-stable subordinator. 
The TSS, denoted by $\{D^{\beta, \nu}_t\}_{t \geq 0}$, with tempering parameter $\nu>0$ and stability index $\beta \in(0,1)$, has the following LT (see e.g. \cite{meerschaert2013transient})
\begin{equation*}
\mathbb{E}\left[e^{-s D^{\beta, \nu}_t}\right]=
e^{-t\big((s + \nu)^{\beta}-\nu^{\beta}\big)}, t\geq 0,s>0. \end{equation*} 
Note that the tempered stable distribution is obtained by exponential tempering of the  $\beta$-stable distribution  (see \cite{Rosinski09}). The advantage of the tempered $\beta$-stable distribution over a $\beta$-stable distribution is that it has moments of all orders and its one-dimensional distributions are infinitely divisible.  The probability density function (pdf) of the one-dimensional distribution of the TSS can be written as (see \cite{meerschaert2013transient})
\begin{equation*}
 f_{\beta, \nu}(x, t)= e^{-\nu x+\nu^{\beta}t} f_{\beta}(x,t),~~ \nu>0, \;\beta\in (0,1), t\geq 0, x\geq 0,
\end{equation*}
where $f_{\beta}(x,t)$ is the marginal pdf of an $\beta$-stable subordinator \cite{UchVla1999}. 
We denote the right-continuous inverse of the TSS, called as the inverse tempered $\beta$-stable subordinator (ITSS), by $\{E^{\beta,\nu}(t)\}_{t \geq 0}$ 
and its pdf by $h_{\beta, \nu}(x,t)$. Various integral representations of the pdf of the ITSS are discussed in \cite{ITS-density, Kumar-ITSS}; for more properties of inverse subordinators see e.g. \cite{LRD2014}.  The LT of the pdf $h_{\beta, \nu}(x, t)$ wrt $`x$' is given by 
 $$
\tilde{h}_{\beta, \nu}(w,t)=\mathcal{L}(h_{\beta, \nu}(x,t): w)=\int_{0}^{\infty}e^{-w x }h_{\beta, \nu}(x,t)dx.
$$
The LT of the above expression wrt `$t$' is given by 
\begin{equation}\label{Item-LT}
 \mathcal{L}(\tilde{h}_{\beta, \nu}(w,t): s)=\frac{f(s)}{s(w +f(s))}.
\end{equation}
\subsection{Generalized fractional derivatives} 
		The generalized Caputo-Djrbashian (C-D) derivative wrt the Bernst\'ein function $f$, which is defined on the space of absolutely continuous functions as follows (see \cite[Definition 2.4]{toaldo2015convolution})
			\begin{equation}\label{Gen_Caputo_FD}
				\frac{\partial^{f}}{\partial t^f}u(t)=b\frac{d}{dt}u(t)+\int \frac{\partial}{\partial t}u(t-s)\Bar{\upsilon}(s)ds,
			\end{equation}
			where $\Bar{\upsilon}(s) = a + \Bar{\upsilon}(s, \infty)$ is the tail of the L\'evy measure and $\upsilon(s,\infty)<\infty,s>0$.\\

If $f(s)=(s+\nu)^\beta-\nu^\beta$ and $\Bar{\upsilon}(s)=\frac{\beta e^{-\nu \beta \Gamma(-\beta, s)}}{\Gamma(1-\beta)}$, where $\Gamma(-\beta, s)=\int_s^{\infty} e^{-z}z^{-\beta-1}dz$, then the generalized fractional derivative  can be called as the C-D tempered fractional derivative of order $\beta\in (0,1)$ with tempering parameter $\nu>0$ is defined by
\begin{equation}\label{CTFD}
\frac{\partial^{\beta,\nu}}{\partial t^{\beta,\nu}}u(t) =  \frac{1}{\Gamma(1-\beta)}\frac{\partial}{\partial t}\int_{0}^{t}\frac{u(r)dr}{(t-r)^{\beta}} - \frac{u(0)}{\Gamma(1-\beta)}\int_{t}^{\infty}e^{-\nu r}\beta r^{-\beta-1}dr.
\end{equation}
The LT for the C-D tempered fractional derivative for a function $u(t)$ satisfies
\begin{equation*}
\mathcal{L}\left[\frac{\partial^{\beta,\nu}}{\partial t^{\beta,\nu}}u\right](s) = ((s+\nu)^{\beta}-\nu^{\beta})\tilde{u}(s) - s^{-1}((s+\nu)^{\beta}-\nu^{\beta}) u(0). 
\end{equation*} 
\subsection{The Hawkes Process} Let $\{N_t\}_{t\geq 0}$ be a counting process with intensity function $\lambda_t$. Now, we present the following definition.
\begin{definition}[The Hawkes Process]
The Hawkes process (HP) (see \cite{dassios2013exact,hainaut2020fractional}) $\{X_t\}_{t\geq 0}=\{(N_t,\lambda_t)\}_{t\geq 0}$ with exponentially decaying intensity is a counting process, defined on $(\Omega,\mathcal{F}_t,\mathbb{P})$, with stochastic intensity  $\{\lambda_t\}_{t\geq 0}$. The intensity is  a self-exciting stochastic process   which depends upon the history of the point process $\{P_t\}_{t \geq 0}$ through the following auto-regressive relation
\begin{align*}
    \lambda_t=\theta+e^{-\kappa(t-s)}(\lambda_s-\theta)+\eta \int_s^t e^{-\kappa(t-u)}dP_u,\;\; t\geq s>0,
\end{align*}
where $\theta\geq 0, \eta>0,$ and $ \kappa>0$ are parameters and  $\{P_t\}_{t\geq0}$ is a continuous-time random walk defined as
\begin{equation*}
P_t = \sum_{i=1}^{N_t} \xi_i, \;\; t\geq 0.
\end{equation*}
Here, $(\xi_i)_{i=1,N_t}$ are sizes of self-excited jumps, a sequence of i.i.d. positive random
variables with finite positive mean $\mu$ and finite variance $\psi^2$.
\end{definition}
The pdf of the HP,  denoted by $p(x,t|y,s)=\mathbb{P}(\lambda_t \in [x, x
+dx]|\lambda_s =y), s\leq t,x>0,y>0$, is the solution of the following forward  differential equation (see \cite{hainaut2020fractional}) 
\begin{align}\label{DDE_HP}
    \frac{\partial p(x,t|y,s)}{\partial t} = &-\frac{\partial}{\partial x} (k(\theta -x)p(x,t|y,s))-\eta \mathbb{E}[\zeta p(x-\eta \zeta,t|y,s)]\nonumber\\&+x\mathbb{E}[p(x-\eta \zeta,t|y,s)-p(x,t|y,s)],
\end{align}
with initial condition $p(x,t|y,s)=\delta_{\{x-y\}},$ where $\delta_x$
is the Dirac measure located at $x$. 
Note that the pdf also solves the backward Kolmogorov equation (see \cite{hainaut2020fractional})
\begin{align}\label{BDE_HP}
    \frac{\partial p(x,t|y,s)}{\partial t} = &k(\theta -y)\frac{\partial p(x,t|y,s)}{\partial y}+y\mathbb{E}[p(x \zeta,t|y+\eta \zeta,s)-p(x,t|y,s)].
\end{align}
The mean $\mathbb{E}_{t}[\cdot]=\mathbb{E}[\cdot|\mathcal{F}_t],t\geq 0$, variance $\mathbb{V}_{t}[\cdot]=\mathbb{V}[\cdot|\mathcal{F}_t],t\geq 0$ and covariance $\mathbb{C}_{t}[\cdot]=\mathbb{C}[\cdot|\mathcal{F}_t],t\geq 0$  are given by  (see \cite{dassios2011dynamic, hainaut2017clustered})
\begin{equation*}
    \mathbb{E}_{0}[\lambda_t] =e^{t(\eta \mu -\kappa)}\lambda_0
+\frac{\kappa \theta}{\eta \mu-\kappa}(e^{t(\eta \mu -\kappa)}-1),\end{equation*}
\begin{equation*}
    \mathbb{V}_{0}[\lambda_t] =\frac{\rho_1\lambda_0+\rho_2}{\eta \mu -\kappa}(e^{2t(\eta \mu -\kappa)}-e^{t(\eta \mu -\kappa)})+\frac{\rho_2}{2(\eta \mu -\kappa)}(1-e^{2(\eta \mu -\kappa)t}),
\end{equation*}
where $\rho_1= \eta^2(\psi^2+\mu^2)$ and $\rho_2= \frac{\eta^2 \kappa \theta(\psi^2+\mu^2)}{\eta \mu -\kappa}$.
The covariance between $\lambda_s$ and $\lambda_t$ for $ s \leq t$ is given by (see \cite{hainaut2020fractional})
$$
\mathbb{C}_0[\lambda_s, \lambda_t]=e^{(\eta \mu -\kappa)(t-s)}\mathbb{V}_0[\lambda_s].
$$ 
The mean and variance are finite when  $\eta \mu-\kappa <0$ (see \cite{hainaut2020fractional} and references therein for more details).

\section{Tempered fractional Hawkes Process}\label{sec:TFHP}

In this section, we define and  discuss some properties of the time-changed Hawkes process.
\begin{definition}[Tempered fractional Hawkes process] The tempered fractional Hawkes process (TFHP) is defined by time-changing the HP with an independent ITSS, and given by
\begin{equation*}
    X({E^{\beta,\nu}_{t}})=\left(N({E_t^{\beta, \nu}}),\lambda({E_t^{\beta, \nu}}) \right), t\geq 0,
\end{equation*}
where $\{E^{\beta,\nu}_t\}_{t\geq0}$ is an independent ITSS.

\end{definition}
\noindent The transition pdf $p^{\beta, \nu}(x,t | y, s)$ of the intensity of jumps is given by
$$
p^{\beta, \nu}(x,t | y, s)= \mathbb{P}(\lambda_{E_t^{\beta, \nu}} \in [x, x+dx]| \lambda_{E_t^{\beta, \nu}}=y).
$$
\begin{theorem}\label{tfhp-pde}
    The pdf $p^{\beta, \nu}(x,t | y, 0)= p^{\beta, \nu}(x,t)$ is the solution of a tempered fractional forward Kolmogorov's equation 
    \begin{align}\label{FDDE_FHP}
     \frac{\partial^{\beta, \nu}p^{\beta, \nu}(x,t | y, 0)}{\partial t^{\beta, \nu}} &=-\frac{\partial}{\partial x}(k (\theta-x)p^{\beta, \nu}(x,t | y, 0))
    -\eta \mathbb{E}[\zeta p^{\beta, \nu}(x-\eta \zeta,t | y, 0)]\nonumber\\
    &+x \mathbb{E}[p^{\beta, \nu}(x-\eta \zeta,t | y, 0)-p^{\beta, \nu}(x,t | y, 0)],  
    \end{align}
    with the condition $p^{\beta, \nu}(x,0 | y, 0)= \delta_{\{x-y\}}$. The pdf also satisfies the tempered fractional backward Kolmogorov's equation
    \begin{align*}
     \frac{\partial^{\beta, \nu} p^{\beta, \nu}(x,t | y, 0)}{\partial t^{\beta, \nu}} &=-k (\theta-y)\frac{\partial p^{\beta, \nu}(x,t | y, 0)} {\partial y}+y \mathbb{E}[p^{\beta, \nu}(x,t | y+\eta \zeta, 0)-p^{\beta, \nu}(x,t | y, 0)],
     \end{align*}
    where $\frac{\partial^{\beta, \nu}}{\partial t^{\beta, \nu}}$ is the C-D tempered fractional derivative of order $\beta \in (0, 1)$, with tempering parameter  $\nu >0$, given in \eqref{CTFD}.
    \end{theorem}
\begin{proof}
     Let $h_{\beta, \nu}(y,t)$ be the pdf of the ITSS  $E_{\beta, \nu}(t)$, and $p(x,t)= p(x, t| y,0) $ be pdf of the HP. Using conditional argument, we have that
		$$
	p^{\beta, \nu}(x,t) = \int_{0}^{\infty} p(x,v)h_{\beta, \nu}(v,t) dv.
			$$
Applying the LT on both sides in the above equation wrt $t$, we get
\begin{align}\label{LT_t_pdf}
\int_{0}^{\infty} p^{\beta, \nu}(x,\gamma) e^{-s t}dt=\Bar{p}^{\beta, \nu}(x, s) &= \int_{0}^{\infty}\int_{0}^{\infty} p(x,u) e ^{-s t}h_{\beta, \nu}(u,t) du dt\nonumber\\
& = \int_{0}^{\infty} p(x,u) \Bar{h}_{\beta, \nu}(u,s)du.
\end{align}
We know that the  LT of $h_{\beta, \nu}(u,t)$ wrt $`t$' is given by 
$\Bar{h}_{\beta, \nu}(v,s) =\int_{0}^{\infty} e ^{-s t}h_{\beta, \nu}(v,t)dt = \frac{f(s)}{s}e^{-v f(s)}, $ 
where $f(s)=(s+\nu)^{\beta}-\nu^{\beta}$. Now, substituting it in equation \eqref{LT_t_pdf} yields
$$
\Bar{p}^{\beta, \nu}(x,w)= \frac{f(s)}{s}\int_{0}^{\infty} p(x,v) e^{-v f(s)} du= \frac{f(s)}{s} \Bar{p}(x,f(s)).
$$
Using the equation \eqref{DDE_HP}, we deduce that $\Bar{p}(x,s)$ is  the solution of 
\begin{align*}
    s\Bar{p}(x, s)-p(x, 0) &=-\frac{\partial}{\partial x}\left(k (\theta-x)\Bar{p}(x, s)\right)-\eta \mathbb{E}[\zeta \Bar{p}(x-\eta \zeta, s)]+x\mathbb{E}[\Bar{p}(x-\eta \zeta, s)-\Bar{p}(x, s)].
\end{align*}
Moreover, replacing $s$ by $f(s)$ leads to
\begin{align*}
     f(s)\Bar{p}(x, f(s))-p(x, 0) &=-\frac{\partial}{\partial x}\left(k (\theta-x)\Bar{p}(x, f(s))\right)-\eta \mathbb{E}[\zeta \Bar{p}(x-\eta \zeta, f(s))]+x\mathbb{E}[\Bar{p}(x-\eta \zeta, f(s))-\Bar{p}(x, f(s))].
\end{align*}
Multiplying by $\frac{f(s)}{s}$ in the above expression, we obtain that 
\begin{align*}
     f(s)\left(\frac{f(s)}{s}\Bar{p}(x, f(s)\right)-&\frac{f(s)}{s}p(x, 0) =-\frac{\partial}{\partial x}\left(k (\theta-x)\frac{f(s)}{s}\Bar{p}(x, f(s))\right)\nonumber\\
    &-\eta \mathbb{E}\left[\zeta \frac{f(s)}{s}\Bar{p}(x-\eta \zeta, f(s))\right]+x\mathbb{E}\left[\frac{f(s)}{s}\Bar{p}(x-\eta \zeta, f(s))-\frac{f(s)}{s}\Bar{p}(x, f(s))\right].
\end{align*}
Since $p^{\beta, \nu}(x,0)= p(x, 0)$, we have that
\begin{align*}
     f(s)\Bar{p}^{\beta, \nu}(x, s)-\frac{f(s)}{s}p^{\beta, \nu}(x, 0) &=-\frac{\partial}{\partial x}\left(k (\theta-x)\Bar{p}^{\beta, \nu}(x, s)\right)-\eta \mathbb{E}\left[\zeta \Bar{p}^{\beta, \nu}(x-\theta \zeta, s)\right]\nonumber\\&+x\mathbb{E}\left[\Bar{p}^{\beta, \nu}(x-\theta \zeta, s)-\Bar{p}^{\beta, \nu}(x, s)\right].
\end{align*}
 It is to note that  the LHS in above equation is the LT of the tempered C-D fractional derivative of $p^{\beta, \nu}(x, t)$. Now, we take the inverse LT on both sides of above equation to get our result in 
\eqref{FDDE_FHP}. \\\\ 
Next, we prove that the given pdf is  also solution of the tempered fractional backward Kolmogorov equation. We begin with rewriting the backward equation \eqref{BDE_HP}
$$
 \frac{\partial p(x,t|y,0)}{\partial t} = k(\theta -y)\frac{\partial p(x,t|y,0)}{\partial y} +y\mathbb{E}[p(x,t|y+\eta \zeta,0)-p(x,t|y,0)].
$$
Using similar approach, we take the LT wrt $t$. Now we follow the same line of argument as above to get the following equation
\begin{align*}
     f(s)\Bar{p}^{\beta, \nu}(x,s|y,0)-&\frac{f(s)}{s}p^{\beta, \nu}(x,0|y,0) =k(\theta -y)\frac{\partial \Bar{p}^{\beta, \nu}(x,s|y,0)}{\partial y}\\&+y\mathbb{E}\left[\Bar{p}^{\beta, \nu}(x,s|y+\eta \zeta,0)-\Bar{p}^{\beta, \nu}(x,s|y,0)\right].
\end{align*}
Hence, the result is proved.
\end{proof}
\noindent  Let ${\varphi}^{\beta,\nu}( )$ be the LT the transition pdf $p^{\beta, \nu}(x, t|y, 0)$ of the time-changed intensity,
    $$
    \varphi^{\beta, \nu}(z, t|y, 0)= \mathbb{E}[e^{-z\lambda_{E^{\beta, \nu}_t}}|\lambda_0=y]= \int_{0}^{\infty} \varphi(z, \tau|y, 0) h_{\beta, \nu}(\tau, t)d \tau,\;\; z>0.
    $$
    \begin{theorem}
    The LT $\varphi^{\beta, \nu}(z, t|y, 0)$ is the solution of the governing equation
    \begin{align*}
    \frac{\partial^{\beta, \nu} \varphi^{\beta, \nu}(z, t|y, 0)}{\partial t^{\beta, \nu}} =-\kappa \gamma z \varphi^{\beta, \nu}(z, t|y, 0)+(1-\kappa z-\mathbb{E}[e^{-z \eta \chi}])\frac{\partial \varphi^{\beta, \nu}(z, t|y, 0)}{\partial z},
\end{align*}
with initial condition $\varphi^{\beta, \nu}(z, 0|y, 0)= e^{-zy}.$
\begin{proof}
    The proof follows a similar approach as in Theorem \eqref{tfhp-pde}.
\end{proof}    \end{theorem}
\noindent Now, we derive some distributional properties for the stochastic intensity function of the TFHP. 
\begin{theorem}\label{mean}
    Let $\gamma=\kappa-\eta\mu>0$, then  mean of $\lambda(E_{t}^{\beta, \nu})$ is given by
\begin{equation*}
    \mathbb{E}_{0}[\lambda(E_{t}^{\beta, \nu})]=\left(\lambda_0-\frac{\kappa \theta}{\gamma}\right) \left(1-\gamma e^{-\nu t}\sum_{m=0}^{\infty}\nu^m t^{\beta+m}M_{\beta, \beta+m+1}^1((\nu^\beta-\gamma) t^{\beta})\right)+\frac{\kappa \theta}{\gamma}.
\end{equation*}
The variance of $\lambda(E_{t}^{\beta, \nu})$ is given by
\begin{align*}
    \mathbb{V}_{0}[\lambda(E_{t}^{\beta, \nu})]=& e^{-\nu t}(\rho_1\lambda_0+\rho_2)\left[\sum_{m=0}^{\infty}2\nu^m t^{\beta+m}M_{\beta, \beta+m+1}((\nu^{\beta}-2\gamma)t^\beta)-\sum_{m=0}^{\infty}\nu^m t^{\beta+m}M_{\beta, \beta+m+1}((\nu^{\beta}-\gamma)t^\beta)\right]\nonumber\\
    &+ \rho_2 e^{-\nu t}\sum_{m=0}^{\infty}\nu^m t^{\beta+m}M_{\beta, \beta+m+1}((\nu^{\beta}-2\gamma)t^\beta)+\left(\lambda_0-\frac{\kappa \theta}{\gamma}\right)^2\mathbb{V}_{0}\left[e ^{-\gamma E_{t}^{\beta, \nu}}\right].
\end{align*}
\begin{proof} Let  $\{\mathcal{H}_t,\}_{t \geq 0}$ be the natural filtration of the ITSS $\{E_{t}^{\beta, \nu}\}_{t\geq0}$.
Using the conditional argument and independence of $\lambda(t)$ and $E^{\beta, \nu}_t$, we have that 
\begin{equation}
    \label{mean-thm1}
\mathbb{E}_0[\lambda(E_{t}^{\beta, \nu})]= \mathbb{E}\left[\mathbb{E}_0[\lambda(E_{t}^{\beta, \nu})|\mathcal{H}_t \vee \mathcal{F}_0)]\right]=\left(\lambda_0-\frac{\kappa \theta}{\gamma}\right)\mathbb{E}[e^{-\gamma E^{\beta, \nu}_t}]+\frac{\kappa \theta}{\gamma}.
\end{equation}

Using \eqref{Item-LT} and taking the inverse LT wrt $`t$', we get
\begin{equation}\label{inver_LT_exp}
\mathbb{E}[e^{-\gamma E^{\beta, \nu}}]= \mathcal{L}^{-1}\left(\frac{f(s)}{s(\gamma+f(s))}: t \right)=\mathcal{L}^{-1} \left(\frac{1}{s}: t\right)-\mathcal{L}^{-1}_t\left(\frac{\gamma}{\gamma+f(s)}:t\right).
\end{equation}
Next, we calculate the inverse LT of  $G(s)$, given by
$$
G(s)=\frac{\gamma}{(\gamma-\nu^\beta) + (s+\nu)^\beta}.
$$
Using the shifting property, we have that 
$$
\mathcal{L}^{-1}(G(s): t)= e^{-\nu t}\mathcal{L}^{-1}\left(\frac{\gamma}{(\gamma-\nu^\beta) + s^\beta}\right).
$$
Now, we apply this result to get (see \cite{kumar2019fractional})
\begin{equation}\label{LT_mitt}
\mathcal{L}^{-1}\left(\frac{s^{ac-b}}{(s^a+\omega)^c}: y\right)= y^{b-1}M_{a, b}^{c}(-\omega y^a), 
\end{equation}
 where $M_{a, b}^{c}(-\omega y^a)$ is generalized Mittag-leffler function (see \cite{prabhakar1971singular}). It is defined by 
  $
 M_{a, b}^{c}(z)=  \sum_{n=0}^{\infty} \frac{(c)_n}{\Gamma(an+b)}\frac{z^n}{n!},
 $
 where $a,b \in \mathbb{C}$ with $\mathcal{R}(a), \;\mathcal{R}(b),\; \mathcal{R}(c) >0$ and $(c)_n=c(c+1)\cdot(c+n-1)$ (for $n=0,1,2,\ldots$ and $c\neq0$) and one otherwise,  is Pochhammer symbol.
Now, inverting the LT of $G(s)$  with the help of \eqref{LT_mitt}, we get
 $$
 \mathcal{L}^{-1}(G(s): t)= \gamma e^{-\nu t} t^{\beta-1}M_{\beta, \beta}^1((\nu^\beta-\gamma) t^{\beta}),
 $$
 which follows by taking $a=\beta, b=\beta, c=1$ and $\omega=(\kappa-\eta \mu-\nu^\beta)$.
 Note that 
 $$
 \mathcal{L}^{-1}\left(\frac{G(s)}{s}: t\right) = \int_{0}^t \gamma e^{-\nu u} u^{\beta-1}M_{\beta, \beta}^1((\nu^\beta-\gamma) u^{\beta}) du,
 $$
We know that (see \cite{kilbas2004generalized})
\begin{equation}\label{kilbas-int}
\int_{0}^{t}y^{d-1}M_{e, d}^{k}(w y^{e})(t-y)^{k-1}dy= \Gamma(k)t^{k+d-1}M_{e, d+k}^{k}(w y^{e}).
\end{equation}
Using the above equation, we have the following 
$$
\mathcal{L}^{-1}\left(\frac{G(s)}{s}: t\right)= \gamma e^{-\nu t}\sum_{m=0}^{\infty}\nu^m t^{\beta+m}M_{\beta, \beta+m+1}^1((\nu^\beta-\gamma) t^{\beta}).
$$
\noindent Now, substitute the above equation in \eqref{inver_LT_exp} to get  \begin{align}\label{meanLT}
\Phi_{\gamma}(t)=\mathbb{E}[e^{-\gamma E^{\beta, \nu}_t}]=1-\gamma e^{-\nu t}\sum_{m=0}^{\infty}\nu^m t^{\beta+m}M_{\beta, \beta+m+1}^1((\nu^\beta-\gamma) t^{\beta}),
\end{align}
and $\Phi_{\gamma}(0)=1$.
The first part of this result can be obtained by substituting the above equation in \eqref{mean-thm1}. 
    Next, we compute the variance of $\lambda(E_{t}^{\beta, \nu})$
    \begin{align*}
\mathbb{V}_0[\lambda(E_{t}^{\beta, \nu})]= &\mathbb{E}\left[\mathbb{V}_0[\lambda(E_{t}^{\beta, \nu})|\mathcal{H}_t \vee \mathcal{F}_0]\right]+ \mathbb{V}_{0}\left[\mathbb{E}[\lambda(E_{t}^{\beta, \nu})|\mathcal{H}_t \vee \mathcal{F}_0]\right]\\
&= \frac{\rho_1\lambda_0+\rho_2}{\gamma}\left[\mathbb{E}\left[e^{-\gamma E_{t}^{\beta, \nu}}\right]-\mathbb{E}\left[e^{-2\gamma E_{t}^{\beta, \nu}}\right]\right]+\frac{\rho_2}{2\gamma}\left[\mathbb{E}\left[e^{-2\gamma E_{t}^{\beta, \nu}}-1\right]\right]\\
& + \left(\lambda_0 -\frac{\kappa \theta}{\gamma}\right)^2\mathbb{V}_{0}\left[e ^{-\gamma E_{t}^{\beta, \nu}}\right].
    \end{align*}
Now, using \eqref{meanLT} and above equation, we complete the proof of the  second part of this Theorem.
\end{proof}
\end{theorem}

\section{Covariance structure of the TFHP}\label{sec:covstruc}
\noindent In this section, we derive the covariance structure of the fractional intensity process of the TFHP. We will first work out some preliminary results in that direction. More specifically, we find out  the analytic expression for the LT of ${E_{s}^{\beta, \nu}+E_{t}^{\beta, \nu}}$ and $E_{t}^{\beta, \nu}-E_{s}^{\beta, \nu}$. It is to note that the correlation structure for the non-tempered is derived in \cite{Leonenko2013,mijena,hainaut2020fractional}, while the following results are not known for the ITSS (its sum/difference) to the best of our knowledge.
\begin{lemma}\label{lemma-1}
    
    For $s\leq t$ and $\gamma>0$. The  LT of $E_{s}^{\beta, \nu}+E_{t}^{\beta, \nu}$ is given by
\begin{align*}
    \mathbb{E}_0\left[e^{-\gamma (E_{t}^{\beta, \nu}+E_{s}^{\beta, \nu})}\right]= -&\gamma \int_{y=0}^{s} h(y)\Phi_\gamma(t-y)dy+e^{-\nu s}\sum_{m=0}^{\infty}\nu^m
s^{\beta +m}M_{\beta, \beta+m+1}(-(2\gamma-\nu^\beta)(s)^{\beta})\nonumber\\
&-\frac{1}{2}+\frac{1}{2}\Phi_{2\gamma}(s)+\Phi_{2\gamma}(t), 
\end{align*}
where $h(y)=e^{-\nu y}y^{\beta-1}M_{\beta,\beta}(\nu^{\beta}y^{\beta})$ and $\Phi_{\gamma}(t)$ is given in \eqref{meanLT}. 
\begin{proof}
       Let $H(u,v)$ be the bivariate cdf of the pair $(E_{s}^{\beta, \nu}, E_{t}^{\beta, \nu})$ respectively. By definition $H(u,\infty)=\mathbb{P}(E_{s}^{\beta, \nu}\leq u)$, $H(\infty,v)=\mathbb{P}(E_{t}^{\beta, \nu}  \leq v)$ and $H(\infty, \infty)=1$. The LT of $E_{t}^{\beta, \nu}+E_{s}^{\beta, \nu}$ is 
       $$
       \mathbb{E}_0\left[e^{-\gamma (E_{t}^{\beta, \nu}+E_{s}^{\beta, \nu})}\right]= \int_{0}^{\infty} \int_{0}^{\infty} e^{-\gamma(u+v)}H(du,dv),
       $$
where 
$H(u,v)= \mathbb{P}(E_{s}^{\beta, \nu} \leq u, E_{t}^{\beta, \nu})\leq v)$. We use the following bivariate integration by parts formula
\begin{align*}
    \int_{0}^{\infty} \int_{0}^{\infty} e^{-\gamma(u+v)}H(du,dv) & = \int_{0}^{\infty} \int_{0}^{\infty}H([u,\infty]\times [v, \infty]f(du,dv)+\int_{0}^{\infty}H([u, \infty]\times[0, \infty])f(du, 0)\\
    &+ \int_{0}^{\infty} H([0, \infty]\times [v,\infty]) f(0, dv)+ f(0,0)H(\infty, \infty),
\end{align*}
with $f(u,v)= e^{-\gamma(u+v)}$. We know that $\mathbb{P}[E_{t}^{\beta, \nu} >0]=1,t>0$, this implies that $f(0,0)H(\infty, \infty)=1$ and $H([u, \infty]\times[0, \infty])= \mathbb{P}(E_{s}^{\beta, \nu}\geq u)$. Therefore, we have that 
\begin{align}\label{bi_ITSS}
    \int_{0}^{\infty} \int_{0}^{\infty} e^{-\gamma(u+v)}H(du,dv) &= \gamma^2 \int_{0}^{\infty} \int_{0}^{\infty}\mathbb{P}(E_{s}^{\beta, \nu} \geq u, E_{t}^{\beta, \nu} \geq v)e^{-\gamma(u+v)}dudv\nonumber\\
&-\gamma \int_{0}^{\infty} e^{-\gamma u}\mathbb{P}(E_{s}^{\beta, \nu}\geq u) du-\gamma \int_{0}^{\infty} e^{-\gamma v}\mathbb{P}(E_{t}^{\beta, \nu}\geq v) dv+1.
\end{align}
The second term of the above equation \eqref{bi_ITSS} can be simplified as follows
\begin{align*}
-\gamma \int_{0}^{\infty} e^{-\gamma u}\mathbb{P}(E_{s}^{\beta, \nu}\geq u) du &= -\gamma \int_{0}^{\infty} e^{-\gamma u}(1-\mathbb{P}(E_{s}^{\beta, \nu}\leq u)du=\mathbb{E}[e^{-\gamma E_{s}^{\beta, \nu})}]+\left[e^{-\gamma u}(1-\mathbb{P}(E_{s}^{\beta, \nu}\leq u)\right]^{\infty}_{0}\\&= \mathbb{E}[e^{-\gamma E_{s}^{\beta, \nu}}]-1.
\end{align*}
Similarly, the third term of \eqref{bi_ITSS} is given by
$$
-\gamma \int_{0}^{\infty} e^{-\gamma v}\mathbb{P}(E_{t}^{\beta, \nu}\geq v) dv=-\gamma \int_{0}^{\infty} e^{-\gamma v}(1-\mathbb{P}(E_{t}^{\beta, \nu}\leq v))dv= \mathbb{E}[e^{-\gamma E_{t}^{\beta, \nu}}]-1.
$$
Hence, equation \eqref{bi_ITSS} can be rewritten as  
\begin{align}\label{final}
 \int_{0}^{\infty} \int_{0}^{\infty} e^{-\gamma(u+v)}H(du,dv)= I+ \mathbb{E}[e^{-\gamma E_{t}^{\beta, \nu}}]+\mathbb{E}[e^{-\gamma E_{s}^{\beta, \nu}}]-1,
\end{align}
 where 
 $$
 I=\gamma^2 \int_{0}^{\infty} \int_{0}^{\infty}\mathbb{P}(E_{s}^{\beta, \nu} \geq u, E_{t}^{\beta, \nu} \geq v)e^{-\gamma(u+v)}dudv.
 $$
We know that the sample paths of  $E_{t}^{\beta, \nu}$ are increasing and continuous. Therefore, for $s \leq t$, we have  $E_{s}^{\beta, \nu} \leq E_{t}^{\beta, \nu}$, and this implies  $\mathbb{P}(E_{t}^{\beta, \nu} \geq v, E_{s}^{\beta, \nu} \geq u) = \mathbb{P}(E_{s}^{\beta, \nu} \geq u)$ for $u >v$. For simplicity of calculations, we write $I=I_1+I_2$, where 
 \begin{align}
     I_1 &= \gamma^2 \int_{0}^{\infty} \int_{0}^{v}\mathbb{P}(E_{s}^{\beta, \nu} \geq u, E_{t}^{\beta, \nu} \geq v)e^{-\gamma(u+v)}dudv, \label{I_1}\\
     I_2 &=\gamma^2 \int_{0}^{\infty} \int_{v}^{\infty}\mathbb{P}(E_{s}^{\beta, \nu} \geq u)e^{-\gamma(u+v)}dudv\nonumber. 
 \end{align}
Rewriting $I_2$ by changing the order of integration, we get
$$
 I_2 = \gamma^2\int_{0}^{\infty} \int_{0}^{u}\mathbb{P}(E_{s}^{\beta, \nu} \geq u)e^{-\gamma(u+v)}dvdu.
 $$
It can be be further simplified as follows
 \begin{align*}
     I_2 &= -\gamma\int_{0}^{\infty}\mathbb{P}(E_{s}^{\beta, \nu} \geq u)(e^{-2u\gamma}-e^{-\gamma u})du\\
     &= \int_{0}^{\infty}(1-\mathbb{P}(E_{s}^{\beta, \nu} \leq u)(-\gamma e^{-2u\gamma})du+ 1-\mathbb{E}[e^{-\gamma E_{s}^{\beta, \nu}})\\
     &= \frac{1}{2}+\frac{1}{2}\mathbb{E}[e^{-2 \gamma E_{s}^{\beta, \nu}}]- \mathbb{E}[e^{-\gamma E_{s}^{\beta, \nu}}].    
 \end{align*}
We know that  $\mathbb{P}(E^{\beta, \nu}_s >u)=\mathbb{P}(D^{\beta, \nu}_u \leq s)$ and   $\{D^{\beta, \nu}_t\}_{t \geq 0}$ has stationary independent increments. Hence, we have that 
 \begin{align}\label{Ind_INCr}
     \mathbb{P}(E^{\beta, \nu}_s >u, E^{\beta, \nu}_t \geq v) &=\mathbb{P}(D^{\beta, \nu}_u \leq s, D^{\beta, \nu}_v \leq t),~~~~ s \leq t\mbox{  and }u \leq v\nonumber\\
    &=  \mathbb{P}(D^{\beta, \nu}_u \leq s, D^{\beta, \nu}_u+ (D^{\beta, \nu}_v- D^{\beta, \nu}_u)\leq t)\nonumber\\
    &=\int_{y=0}^s f_{\beta, \nu}(y, u) \int_{0}^{t-y}f_{\beta, \nu}(x, v-u)dxdy.
 \end{align}
 Substituting the above expression into \eqref{I_1} and using Fubini's theorem, it follows that
 \begin{align}
     I_1 &= \gamma^2 \int_{y=0}^{s}\int_{x=0}^{t-y}\int_{u=0}^{\infty}f_{\beta, \nu}(y, u)\int_{v=u}^{\infty}f_{\beta, \nu}(x, v-u)e^{-\gamma(u+v)} dvdudxdy\nonumber\\
    &= \gamma^2 \int_{y=0}^{s}\int_{x=0}^{t-y}\int_{u=0}^{\infty}f_{\beta, \nu}(y, u) e^{-2\gamma u}du\int_{\gamma=0}^{\infty}f_{\beta, \nu}(x, w)e^{-w\gamma} dwdxdy. \label{I_1form}
 \end{align}
To simplify above calculations, let  $h(y)= \int_{u=0}^{\infty}f_{\beta, \nu}(y, u) e^{-2\gamma u}du$ and $g(x, \gamma)= \int_{\gamma=0}^{\infty}f_{\beta, \nu}(x, w)e^{-w\gamma} dw$. Then, the LT of $h(y)$ is given by
 \begin{align*}
     \mathcal{L}(h(y):s) &= \int_{y=0}^{\infty}e^{-sy}\int_{u=0}^{\infty}f_{\beta, \nu}(y, u) e^{-2\gamma u}dudy\\
     &= \int_{u=0}^{\infty}e^{-2\gamma u}\int_{y=0}^{\infty}e^{-sy}f_{\beta, \nu}(y, u) dydu\\
&=\int_{u=0}^{\infty}e^{-u (2\gamma+(s+\nu)^\beta-\nu^{\beta}})dy\\
&= \frac{1}{(s+v)^\beta+(2\gamma-\nu^\beta)}.
 \end{align*} 
 Now, inverting the LT of the  $h(y)$ by using the equation \eqref{LT_mitt}, is given by
 \begin{equation}\label{LT_h}
 h(y)= e^{-\nu y} y^{\beta-1}M_{\beta, \beta}^1(-(2\gamma-\nu^{\beta})y^{\beta}) 
 \end{equation}
 Next, we find the LT of the $g(x, \gamma)$.
 \begin{align}
    \mathcal{L}(g(x, \gamma): s) &= \int_{w=0}^{\infty} e^{-w\gamma}\int_{x=0}^{\infty}e^{-sx}f_{\beta,\nu}(x, w)dxdw\nonumber\\
    &= \frac{1}{\gamma +((s+\nu)^\beta-\nu^{\beta})}\nonumber\\
    &= \frac{1}{\gamma}\left(1-\frac{(s+\nu)^\beta-\nu^{\beta}}{\gamma +((s+\nu)^\beta-\nu^{\beta})}\right)\nonumber\\
    &= \mathcal{L}\left(\frac{-1}{\gamma}\frac{d}{dx}\Phi_{\gamma}(x): s\right).\label{LT_g}
 \end{align}
 Substituting \eqref{LT_h} and \eqref{LT_g} in \eqref{I_1form}, we have
\begin{align*}
    I_1 &= \gamma^2 \int_{y=0}^{s}\int_{x=0}^{t-y}h(y) \left(\frac{-1}{\gamma}\frac{d}{dx}\Phi_\gamma(x)\right) dxdy\\
    &= -\gamma \int_{y=0}^{s} h(y)(\Phi_\gamma(t-y)-1)dy\\
    &= -\gamma \int_{y=0}^{s} h(y)\Phi_\gamma(t-y)dy +\gamma \int_{y=0}^{s}h(y) dy.
\end{align*}
Recall from equation \eqref{kilbas-int}

$$
\int_{0}^{s} h(y) dy = e^{-\nu s}\sum_{m=0}^{\infty}\nu^m
s^{\beta +m}M_{\beta, \beta+m+1}(-(2\gamma-\nu^\beta)(s)^{\beta}).$$ Then 
$$
I_1= -\gamma \int_{y=0}^{s} h(y)\Phi_\gamma(t-y)dy+e^{-\nu s}\sum_{m=0}^{\infty}\nu^m
s^{\beta r+m}M_{\beta, \beta+m+1}(-(2\gamma-\nu^\beta)s^{\beta}).
$$
 Now, we have 
 \begin{align*}
     I=-&\gamma \int_{y=0}^{s} h(y)\Phi_\gamma(t-y)dy+e^{-\nu s}\sum_{m=0}^{\infty}\nu^m
s^{\beta r+m}M_{\beta, \beta+m+1}(-(2\gamma-\nu^\beta)(s)^{\beta})\nonumber\\
&+\frac{1}{2}+\frac{1}{2}\mathbb{E}[e^{-2 \gamma E_{s}^{\beta, \nu}}) - \mathbb{E}[e^{-\gamma E_{s}^{\beta, \nu}}].
\end{align*}
Substitute the above expression in equation \eqref{final}, we get the desired result.
\end{proof}
\end{lemma}
We now derive another result which will be useful in finding out covariance of the TFHP.
\begin{lemma}\label{lemma-2}
  Let $t\geq s>0$ and $\gamma>0$.  The LT of the bivariate distribution function $H(u, v)= \mathbb{P}[E_t^{\beta, \nu} \leq u, E_s^{\beta, \nu} \leq v ]$ of the process $\{E_t^{\beta, \nu}\}_{t\geq0}$ is given by
    \begin{align*}
         \mathbb{E}_0\left[e^{-\gamma (E_{t}^{\beta, \nu}-E_{s}^{\beta, \nu})}\right]=& \gamma \int_0^{s}h(y)\Phi_{\gamma}(t-y)dy-\gamma e^{-\nu s}\sum_{m=0}^{\infty}\nu^m
s^{\beta +m}M_{\beta, \beta+m+1}(\nu^\beta s^{\beta})\\&+\Phi_{\gamma}(t)+\gamma \mathbb{E}[E^{\beta, \nu}_s], \nonumber
    \end{align*}
 where $h(y)=e^{-\nu y}y^{\beta-1}M_{\beta, \beta}(\nu^{\beta}t^{\beta})$ and $\Phi_{\gamma}(t)$ is  as given  in \eqref{Item-LT}.
\end{lemma}
    \begin{proof}
    We start with considering the following LT of the difference of ITSS
    
       $$
       \mathbb{E}_0\left[e^{-\gamma (E_{t}^{\beta, \nu}-E_{s}^{\beta, \nu})}\right]= \int_{0}^{\infty} \int_{0}^{\infty} e^{-\gamma|u-v|}H(du,dv),
       $$
Analogues to the  computation mentioned  in \cite{Leonenko2013}, we have that
\begin{align}\label{Com_Exp}
    \int_{0}^{\infty} \int_{0}^{\infty} e^{-\gamma|u-v|}H(du, dv) & = I+\gamma \mathbb{E}[E^{\beta, \nu}_s]+\Phi_{\gamma}(t),
\end{align}
where
\begin{align*}
    I=-\gamma^2\int_{u=0}^{\infty}\mathbb{P}[E^{\beta, \nu}_t \geq u, E^{\beta, \nu}_s \geq v]\int_{v=u}^{\infty}e^{-\gamma(u-v)}dvdu.
\end{align*}
It is important to note that the ITSS $\{E^{\beta, \nu}_t\}_{t\geq 0}$ does not have a self-similar property and therefore we take the LT approach to further simplify the above integral. We now substitute equation \eqref{Ind_INCr} into the above expression and applying Fubini's theorem to get
\begin{align}\label{I_form}
     I = -\gamma^2 \int_{y=0}^{s}\int_{x=0}^{t-y}\int_{u=0}^{\infty}f_{\beta, \nu}(y, u)du\int_{\gamma=0}^{\infty}f_{\beta, \nu}(x, w)e^{-\gamma\gamma} dw dx dy.
 \end{align}
  Let $h(y)= \int_{u=0}^{\infty}f_{\beta, \nu}(y, u) du$ and $g(x, \gamma)= \int_{\gamma=0}^{\infty}f_{\beta, \nu}(x, w)e^{-w\gamma} dw$. The LT of $h(y)$ is given by
 \begin{align*}
     \mathcal{L}(h(y):s) &= \int_{y=0}^{\infty}e^{-sy}\int_{u=0}^{\infty}f_{\beta, \nu}(y, u)dudy\nonumber\\
     &= \int_{u=0}^{\infty}\int_{y=0}^{\infty}e^{-sy}f_{\beta, \nu}(y, u) dydu\nonumber\\
&= \frac{1}{(s+\nu)^\beta-\nu^\beta}.
 \end{align*}
The inverse LT of the above equation from equation \eqref{LT_mitt}, we have 
 \begin{equation}\label{LT_h_I}
 h(y)= e^{-\nu y} y^{\beta-1}M_{\beta, \beta}^1(\nu^{\beta}y^{\beta}).
 \end{equation}
 Further, we find the LT of the $g(x, \gamma)$, which has already obtained in equation \eqref{LT_g}, we get
 \begin{align}\label{LT_G_I}
    \mathcal{L}(g(x, \gamma): s) = \mathcal{L}\left(\frac{-1}{\gamma}\frac{d}{dx}\Phi_\gamma(x): s\right).
 \end{align}
 Substituting \eqref{LT_h_I} and \eqref{LT_G_I} in \eqref{I_form}, 
  we obtain the following
\begin{align*}
    I &= -\gamma^2 \int_{y=0}^{s}\int_{x=0}^{t-y}h(y) \left(\frac{-1}{\gamma}\frac{d}{dx}\Phi_\gamma(x)\right) dxdy\\
    &= \gamma \int_{y=0}^{s} h(y)(\Phi_\gamma(t-y)-1)dy\\
    &= \gamma \int_{y=0}^{s} h(y)\Phi_\gamma(t-y)dy -\gamma \int_{y=0}^{s}h(y) dy.
\end{align*}
Recall from equation \eqref{kilbas-int}
$$
\int_{0}^{s} h(y) dy = e^{-\nu s}\sum_{m=0}^{\infty}\nu^m
s^{\beta +m}M_{\beta, \beta+m+1}(\nu^\beta s^{\beta}).$$ Then 
$$
I= \gamma \int_{y=0}^{s} h(y)\Phi_\gamma(t-y)dy-e^{-\nu s}\sum_{m=0}^{\infty}\nu^m
s^{\beta +m}M_{\beta, \beta+m+1}(\nu^\beta s^{\beta}).
$$
Now the result  follows using \eqref{Com_Exp} and the above equation, which completes the proof.
\end{proof}
\noindent We now derive the main result of this section by finding out the covariance  of the stochastic intensity of the TFHP.

    \begin{theorem}\label{cov}
 Let $\gamma>0$, the covariance between $\lambda(E_{t}^{\beta, \nu})$ and $\lambda(E_{s}^{\beta, \nu})$ is given by 
    \begin{align*}
        \mathbb{C}_0(\lambda(E_{t}^{\beta, \nu}), \lambda(E_{s}^{\beta, \nu}))&= 
        \frac{\rho_1\lambda_0+\rho_2}{\eta \mu -\kappa}\left[\mathbb{E}_0\left[e^{-\gamma (E_{t}^{\beta, \nu}+E_{s}^{\beta, \nu})}\right]-\Phi_\gamma(t)\right]\nonumber\\
       & +\frac{\rho_2}{2(\eta \mu -\kappa)}\left[\mathbb{E}_0\left[e^{-\gamma (E_{t}^{\beta, \nu}-E_{s}^{\beta,\nu} )}\right]-\mathbb{E}_0\left[e^{-\gamma (E_{t}^{\beta, \nu}+E_{s}^{\beta, \nu})}\right]\right]\nonumber\\
& + \left(\lambda_0 +\frac{\kappa \theta}{\eta \mu-\kappa}\right)^2\left[\mathbb{E}_0\left[e^{-\gamma (E_{t}^{\beta, \nu}+E_{s}^{\beta, \nu})}\right]- \Phi_{\gamma}(t)\Phi_{\gamma}(s)\right].
    \end{align*}
\begin{proof}
We will provide a brief sketch of the proof and omit the details since they are analogues the proof mentioned in \cite[Proposition 8.2]{hainaut2020fractional}. Let $s\leq t$, consider the following
\begin{align}\label{cov-eq}
\mathbb{C}_0[\lambda(E_{s}^{\beta, \nu}), \lambda(E_{t}^{\beta, \nu})]&= \mathbb{E}\left[\mathbb{C}_0[\lambda(E_{s}^{\beta, \nu}), \lambda(E_{t}^{\beta, \nu})]|\mathcal{H}_t \vee \mathcal{F}_0)\right]\nonumber\\&+ \mathbb{C}_{0}\left[\mathbb{E}[\lambda(E_{s}^{\beta, \nu})|\mathcal{H}_t \vee \mathcal{F}_0], \mathbb{E}[\lambda(E_{t}^{\beta, \nu}), |\mathcal{H}_t \vee \mathcal{F}_0]\right].
\end{align}
\noindent Here,
\begin{align*}
   \mathbb{E}\left[\mathbb{C}_0[\lambda(E_{s}^{\beta, \nu}), \lambda(E_{t}^{\beta, \nu})]|\mathcal{H}_t \vee \mathcal{F}_0)\right] = &\frac{\rho_1\lambda_0+\rho_2}{\eta \mu -\kappa}\left[\mathbb{E}_0[e^{-\gamma (E_{t}^{\beta, \nu}+E_{s}^{\beta, \nu})}]-\mathbb{E}_0[e^{-\gamma (E_{s}^{\beta, \nu})}]\right]\\
   &+\frac{\rho_2}{2(\eta \mu -\kappa)}\left[\mathbb{E}_0\left[e^{-\gamma (E_{t}^{\beta, \nu}-E_{s}^{\beta,\nu} )}\right]-\mathbb{E}_0\left[e^{-\gamma (E_{t}^{\beta, \nu}+E_{s}^{\beta, \nu})}\right]\right].
\end{align*}
Now, using Lemma  \ref{lemma-1} and \ref{lemma-2}, we obtain the first term of the LHS of equation \eqref{cov-eq}. Considering the second term and following similar steps (replacing the inverse $\beta$-stable subordinator with the ITSS) as mentioned in \cite[Proposition 8.2]{hainaut2020fractional}, we obtain the desired result. 
\end{proof}    
\end{theorem}

\section{Generalized fractional Hawkes process}\label{sec:GFHP}
In previous sections, we have considered the ITSS as a time-change in the stochastic intensity of the HP to define the TFHP. However, one can easily define generalized version of the ``fractional" Hawkes process  by taking a general inverse L\'evy process $\{E_f(t)\}_{t\geq0}$, defined in \eqref{inverse-sub}, as time-change for the HP. In this  way, the results obtained in this section applies to a very large class of time-changed Hawkes processes, for example, this definition encapsulates the HP time-changed by inverse of inverse Gaussian subordinator (see \cite{Kumar-Hitting}), inverse gamma (see \cite{invgamma}), inverse stable (see \cite{mnv}), inverse tempered stable (see \cite{Kumar-ITSS,itss-gupta,ITS-density}),  inverse Dickman (see \cite{Gupta-dickman}) and several mixtures of inverse subordinators, \textit{etc}.  

We now define the generalized fractional Hawkes process (GFHP) by using the inverse L\'evy subordinator as a time-change in the HP $\{X(t)\}_{t\geq 0}$

\begin{definition}[Generalized fractional Hawkes process] The generalized fractional Hawkes process (GFHP) is defined by time-changing the HP with an independent inverse Le\'vy subordinator (defined in \eqref{inverse-sub}) and is given by
\begin{equation*}
    X({E_f(t)})=\left(N({E_f(t)}),\lambda({E_f(t)}) \right),t\geq 0.
\end{equation*}

\end{definition}
\noindent We now present some results for the GFHP. \\
The transition pdf $p^{f}(x,t | y, s),s>0,t>0,$ of the stochastic  intensity can be defined as 
$$
p^{f}(x,t | y, s)= \mathbb{P}(\lambda_{E_f(t)} \in [x, x+dx]| \lambda_{E_f(t)}=y),s\leq t,x,y>0.
$$
We now present the governing  difference-differential equation for the one-dimensional distributions of the GFHP.
\begin{theorem}
    The \textit{pdf} $p^{f}(x,t | y, 0):= p^{f}(x,t)$ is the solution of the following  fractional forward difference--differential equation 
    \begin{align*}\label{GFDDE_FHP}
     \frac{\partial^{f}p^{f}(x,t | y, 0)}{\partial t^{f}} &=-\frac{\partial}{\partial x}(k (\theta-x)p^{f}(x,t | y, 0))
    -\eta \mathbb{E}[\zeta p^{f}(x-\eta \zeta,t | y, 0)]\nonumber\\
    &+x \mathbb{E}[p^{f}(x-\eta \zeta,t | y, 0)-p^{f}(x,t | y, 0)],  
    \end{align*}
    with the condition $p^{f}(x,0 | y, 0)= \delta_{\{x-y\}}$. The \textit{pdf} also solves the generalized C-D fractional backward Kolmogorov's equation:
    \begin{align*}
     \frac{\partial^{f} p^{f}(x,t | y, 0)}{\partial t^{f}} &=-k (\theta-y)\frac{\partial p^{f}(x,t | y, 0)} {\partial y}+y \mathbb{E}[p^{f}(x,t | y+\eta \zeta, 0)-p^{f}(x,t | y, 0)],
     \end{align*}
    where $\frac{\partial^{f}}{\partial t^{f}}$ is the generalized C-D derivative, given in \eqref{Gen_Caputo_FD}.
    \end{theorem}
\begin{proof}

We know that the LT of the pdf $h_{f}(x, t)$ of the one-dimensional distribution of $\{E_t^f\}_{t \geq 0}$ wrt $`t$', and  is given by 
 $$ \mathcal{L}(h_{f}(x,t): s)=\frac{f(s)}{s}e^{-xf(s)}.
$$
Consider the above-mentioned LT and retracing the proof of Theorem \ref{tfhp-pde} to get the desired result.
\end{proof}
\noindent We next present mean, variance and covariance of the GFHP.
\begin{theorem}
    The  mean of $\lambda(E_f(t))$ is given
\begin{equation*}
    \mathbb{E}_{0}[\lambda(E_f(t))]=\left(\lambda_0-\frac{\kappa \theta}{\gamma}\right)\mathbb{E}[e^{-\gamma E_f(t)}]+\frac{\kappa \theta}{\gamma}.
\end{equation*}
The variance of $\lambda(E_f(t))$ is given
\begin{align*}
    \mathbb{V}_{0}[\lambda(E_f(t))]=&\frac{\rho_1\lambda_0+\rho_2}{\gamma}\left[\mathbb{E}[e^{-\gamma E_f(t)}]-\mathbb{E}[e^{-2\gamma E_f(t)}]\right]+\frac{\rho_2}{2\gamma}\left[\mathbb{E}[e^{-2\gamma E_f(t)}-1]\right] + \left(\lambda_0 -\frac{\kappa \theta}{\gamma}\right)^2\mathbb{V}_{0}\left[e^{-\gamma E_f(t)}\right].
\end{align*}
The covariance between $\lambda(E_f(t))$ and $\lambda(E_f(s))$ is given by 
    \begin{align*}
        \mathbb{C}_0(\lambda(E_f(t)), \lambda(E_f(s)))&= 
        \frac{\rho_1\lambda_0+\rho_2}{\eta \mu -\kappa}\left[\mathbb{E}_0\left[e^{-\gamma (E_f(t)+E_f(s))}\right]-\mathbb{E}_0\left[e^{-\gamma E_f(t)}\right]\right]\nonumber\\
       & +\frac{\rho_2}{2(\eta \mu -\kappa)}\left[\mathbb{E}_0\left[e^{-\gamma (E_f(t)-E_f(s) )}\right]-\mathbb{E}_0\left[e^{-\gamma (E_f(t)+E_f(s))}\right]\right]\nonumber\\
& + \left(\lambda_0 +\frac{\kappa \theta}{\eta \mu-\kappa}\right)^2\left[\mathbb{E}_0\left[e^{-\gamma (E_f(t)+E_f(s))}\right]- \mathbb{E}_0\left[e^{-\gamma E_f(t)}\right]\mathbb{E}_0\left[e^{-\gamma E_f(s)}\right]\right].
    \end{align*}
\end{theorem}
\begin{proof}
    The proof runs on  similar lines as Theorems \ref{mean}
 and \ref{cov} by replacing the ITSS with the inverse L\'evy subordinator at appropriate places, and hence omitted here. \end{proof}

\def\cprime{$'$}

\end{document}